\documentclass{amsart}
\usepackage[cp1250]{inputenc}

\newtheorem{theorem}{Theorem}[section]

\newtheorem{example}[theorem]{Example}
\newtheorem{remark}[theorem]{Remark}
\newtheorem{conjecture}[theorem]{Conjecture}
\numberwithin{equation}{section}

\author{S. Capparelli}
\address[Stefano Capparelli]
{Dipartimento SBAI,
Universit\`a di Roma La Sapienza,
Roma, Italy}
\email{stefano.capparelli@uniroma1.it}

\author{A. Meurman}
\address[Arne Meurman]{Department of Mathematics, University of Lund, Box 118, 22100 Lund, Sweden}
\email{arne.meurman@math.lu.se}

\author{A. Primc}
\address[Andrej Primc]{Kersnikova 11, 1\,000 Ljubljana, Slovenia}
\email{aprimc@gmail.com}

\author{M. Primc}
\address[Mirko Primc]{Faculty of Science,  University of Zagreb,  Zagreb, Croatia}
\email{primc@math.hr}

\begin{document}

\title
 {New partition identities from $C^{(1)}_\ell$-modules}

\begin{abstract}
In this paper we conjecture combinatorial Rogers-Ramanujan type colored partition identities related to standard representations of the affine Lie algebra of type $C^{(1)}_\ell$, $\ell\geq2$, and  we conjecture  similar colored partition identities with no obvious connection to representation theory of affine Lie algebras.
\end{abstract}

 \maketitle

  \section{Introduction}

\noindent
The seminal work \cite{LW} of J. Lepowsky and R. Wilson on a Lie-theoretic interpretation of the Rogers-Ramanujan identities led to the discovery of numerous new combinatorial identities, like in \cite{C} or \cite{MP}. Recently several identities in the style of the Rogers-Ramanujan identities related to the representation theory of affine Lie algebras have appeared, let us mention only \cite{KR} and \cite{BJSM}. On the other side, some parts of representation theory lead to Rogers-Ramanujan type {\em colored} partition identities, let us mention
 only \cite{DK} and \cite{PS1} for this vein of research.

In this paper we conjecture combinatorial Rogers-Ramanujan type colored partition identities related to standard representations of the affine Lie algebra of type $C^{(1)}_\ell$, $\ell\geq2$, and  we conjecture  similar colored partition identities with no obvious connection to representation theory of affine Lie algebras.

In Section  \ref{Odd array width w} we start with the array of natural numbers $\mathcal N_5$ composed as a multiset of two copies of the set of natural numbers $\mathbb N$ and the additional set of odd numbers, arranged in $5$ rows, with diagonals of width $w=5$. Such an array appears naturally in the representation theory of the affine Kac-Moody Lie algebra $\widehat{\mathfrak{g}}$ of type $C^{(1)}_2$  and it is expected that colored partitions on $\mathcal N_5$, satisfying certain difference and initial conditions, should parametrize bases of standard $\widehat{\mathfrak{g}}$-modules---see Remark \ref{R:known w odd cases} below. In Conjecture \ref{C:conjecture} we guess, inspired by \cite{T}, a possible form of these colored partitions and, by a computing experiment, we conjecture the corresponding Rogers-Ramanujan type colored partition identities. Numerical evidence supports the conjecture for colored partitions on the array $\mathcal N_{2\ell+1}$, related to all standard modules for affine Lie algebra of type $C^{(1)}_\ell$, $\ell\geq2$.

In Section  \ref{Even array width w} we start with the array of natural numbers $\mathcal N_4$ composed as a multiset of two copies of the set of natural numbers $\mathbb N$, arranged in $4$ rows, with diagonals of width $w=4$. In analogy with the ``$w$ odd case'', in Conjecture \ref{C:conjecture w even} we conjecture similar colored partition identities for all $w=2\ell$, $\ell\geq2$, with no obvious connection to representation theory of affine Lie algebras, but again supported by numerical evidence.

The arrays of natural numbers
$$
\mathcal N_2,\quad \mathcal N_3,\quad \mathcal N_4,\quad \mathcal N_5,\quad \mathcal N_6,\quad \mathcal N_7,\quad\dots
$$
form a natural sequence and for each  $\mathcal N_{2\ell}$ and  $\mathcal N_{2\ell+1}$, $\ell\geq1$, and all nonnegative integers $k_0, k_1, \dots, k_\ell$, $k=k_0+ k_1+ \dots+ k_\ell>0$, we have a class of
$(k_0, k_1, \dots, k_\ell)$-admissible colored partitions (defined in (\ref{E:array of frequencies})--(\ref{E:difference condition}) and (\ref{E:array of frequencies on w even})--(\ref{E:difference condition on w even}) below) for which we conjecture Rogers-Ramanujan type combinatorial identities. For $w=2\ell=2$ and $k=1$ we have the two Rogers-Ramanujan identities, and for $w=2\ell+1=3$ and $k=1$ we have two identities (in some sense) equivalent to the two Capparelli identities\footnote{See Remarks \ref{R:known w odd cases}  and \ref{R: Rogers-Ramnujan identities} below.}.

In Section \ref{An algorithm} we describe an algorithm for constructing admissible colored partitions, and in the Appendix we give a Python code for counting admissible colored partitions.

\section{Lepowsky's product formula}

In this section we give the Lie theoretic origin of the product expression in Conjecture \ref{C:conjecture} below. Note however that Conjecture \ref{C:conjecture} is a purely combinatorial statement.

By $C_{\ell}^{(1)}$ we denote the affine Lie algebra of type $C_{\ell}^{(1)}$, see e.g. \cite{K}, Table Aff 1, page 54 and (7.2.1), page 98. By $L(k_0,\dots,k_\ell)$ we denote the irreducible highest weight module $L(\Lambda)$ with highest weight $\Lambda = k_0\Lambda_0+\cdots +k_\ell \Lambda_\ell$, see e.g. \cite{K}, page 147. The principally specialized character $\text{ch}_q L_{C_\ell^{(1)}}(k_0,\dots,k_\ell)$ is defined in \cite{K}, page 152 together with \S10.9, page 181. Associated to $L(k_0,\dots,k_\ell)$ set $k = k_0 + k_1 + \cdots + k_\ell$. 

Let $s_0,s_1,\dots,s_\ell \in\mathbb{N}, s = s_0 + s_1 + \cdots + s_\ell$. We define the following triangular scheme of natural numbers:
\bigskip

\noindent
$D(s_0,s_1, \dots, s_\ell)=\{s_0, \dots, 2s-s_0 \}$ has $2\ell+1$ ascending numbers with increments
\newline$ \qquad s_1,s_2,\dots,s_\ell,s_\ell,\dots,s_2,s_1$,

\noindent
$D(s_1, \dots, s_\ell)=\{s_1, \dots \}$ has $2\ell-1$ ascending numbers with increments
\newline$s_2,\dots,s_\ell,s_\ell,\dots,s_2$,

\noindent
$D(s_2, \dots, s_\ell)=\{s_2, \dots \}$ has $2\ell-3$ ascending numbers with increments
\newline$s_3,\dots,s_\ell,s_\ell,\dots,s_3$,

$\vdots$

\noindent
$D(s_{\ell-1}, s_\ell)=\{s_{\ell-1}, s_{\ell-1}+s_\ell, s_{\ell-1}+2s_\ell \}$,

\noindent
$D(s_\ell)=\{s_\ell \}$.
\bigskip

\begin{example}\label{Ex:congruence triangle (3,2,1,1,2)}
For example, for $(s_0, s_1, s_2, s_3, s_4)=(3,2,1,1,2)$ we have
$$
\begin{matrix}
D(3,2,1,1,2)\\
D(2,1,1,2)\\
D(1,1,2)\\
D(1,2)\\
D(2)\\
\end{matrix}\qquad=\qquad
\begin{matrix}
3&5&6&7&9&11&12&13&15\\
&2&3&4&6&8&9&10&\\
&&1&2&4&6&7&&\\
&&&1&3&5&&&\\
&&&&2&&&&\\
\end{matrix}.
$$
\end{example}

\noindent
Note that we can obtain $D(s_{r+1}, \dots, s_\ell)$ from $D(s_{r}, \dots, s_\ell)$ by eliminating endpoints in  $D(s_{r}, \dots, s_\ell)$ and then subtracting $s_{r}$ from the remaining elements. We define {\it the congruence triangle} as the multiset
$$
\Delta(s_1, \dots, s_\ell)=D(s_1, \dots, s_\ell)\cup D(s_2, \dots, s_\ell)\cup\dots\cup D(s_\ell)
$$
and we denote by $\{0\}^\ell$ the multiset consisting of $\ell$ copies of $0$.

Now we can state {\it Lepowsky's product formula} for $C\sp{(1)}_\ell$, $\ell\geq2$ (cf.  \cite{L}, \cite{B}):
\begin{equation}\label{E:Lepowsky's product formula}
\begin{aligned}\notag
&\text{ch}_qL_{C\sp{(1)}_\ell}(k_0, k_1, \dots, k_\ell)=\\
&\frac{\prod_{a\in\{0\}^\ell\cup D(k_0+1, k_1+1,\dots, k_\ell+1);\,
b\in \Delta( k_1+1,\dots, k_\ell+1);\,
j\equiv a,\pm b\mod (2\ell+2k+2)}(1-q^j)
}
{\prod_{j\, \text{odd}}(1-q^j)\prod_{j\in\mathbb N}(1-q^j)^\ell}.
\end{aligned}
\end{equation}
\begin{example}
By using Example \ref{Ex:congruence triangle (3,2,1,1,2)} we see that
the principally specialized character $\text{ch}_qL_{C\sp{(1)}_4}(2,1,0,0,1)$ is equal to
$$
\prod_{j\, \text{odd}; j\equiv 1,1,2,3,4,4,5,5,6,7,7,8,8,9,10,10,11,11,12,13,13,14,14,15,16,17,17 \mod 18}(1-q^j)^{-1}.
$$
Note that in this product the factor $(1-q^1)^{-1}$ appears three times.
\end{example}

\section{Arrays with odd width  $w\geq5$}\label{Odd array width w}

Let $\mathcal N=\mathcal N_5=\mathcal N_{C^{(1)}_2}$ be the array of natural numbers
\begin{equation}\label{E:C(1)2 principally specialized array}
\begin{matrix}
1& & 3& & 5& & 7& &  {9}&  \\
& 2& &4 & &6 & & {8} & & 10 \\
1& &3 & &5 & &7 & &{9} & \\
&2 & &4 & &6 & &{8} & &10  \\
1& &3 & &5 & &{7} & & {9} &
\end{matrix}\quad\dots \,.
\end{equation}
This array consists of two copies of the set of natural numbers $\mathbb N$ and the additional set of odd numbers, arranged in $5$ rows, with diagonals of width $w=5$. Numbers increase by one going to the right on any diagonal. We shall consider elements in these sets as different, say ``colored'' by their position in the array. For example, the number $7$ appears three times on three different places of the array $\mathcal N$ and we consider $7$ in the first row different from the other two. We say that two elements in an array are {\it adjacent} if they are simultaneously on two adjacent rows and two adjacent diagonals. For example, $6$ and $8$ in the second row are adjacent to $7$ in the first row and, just as well, adjacent to $7$ in the third row. We say that the set $\{a_1,a_2,a_3,\dots\}$ is a downward path $\mathcal Z$ in an array if $a_i$ is in the $i$-th row and if $(a_i, a_{i+1})$ is a pair of two adjacent elements for all $i$. For example, $\mathcal Z=\{7,6,5,4,5\}$ is a downward path in $\mathcal N$ and there are altogether $2^4$ downward paths through $7$ in the first row. In Section \ref{An algorithm} we shall consider downward paths which start from the top row, but need not reach the bottom row. $\mathcal Z=\{7,6,5,4\}$ is an example of such downward path---we shall say that this $\mathcal Z=\{7,6,5,4\}$ ends at $4$ in the fourth row.

We consider colored partitions
\begin{equation}\label{E:colored partition on the array}
n=\sum_{a\in\mathcal N}f_a\cdot a,
\end{equation}
where $f_a$ is the frequency of the part $a\in\mathcal N$ in the colored partition (\ref{E:colored partition on the array}) of $n$.
Let $k_0, k_1, k_2\in\mathbb N_0$, $k=k_0+k_1+k_2>0$. We say that an array of frequencies    $\mathcal F$
\begin{equation}\label{E:array of frequencies}
\qquad\qquad
\begin{matrix}
f_1& &f_ 3& & f_5& & f_7&   \\
&f_ 2& &f_4 & &f_6 & & f_{8}  \\
f_1& &f_ 3& & f_5& & f_7&   \\
&f_ 2& &f_4 & &f_6 & & f_{8}  \\
f_1& &f_ 3& & f_5& & f_7&    \\
\end{matrix}\quad\dots
\end{equation}
{\it is $(k_0, k_1, k_2)$-admissible} if {\it the extended array of frequencies} $\mathcal F\sp{(k_0, k_1, k_2)}$
\begin{equation}\label{E:extended array of frequencies}
\begin{matrix}
 k_{2}&& f_1& &f_ 3& & f_5& & f_7&    \\
 & 0&&f_ 2& &f_4 & &f_6 & & f_{8}  \\
 k_{1}&& f_1& &f_ 3& & f_5& & f_7&   \\
 & 0&&f_ 2& &f_4 & &f_6 & & f_{8}  \\
 k_{0}&& f_1& &f_ 3& & f_5& & f_7&    \\
\end{matrix}\quad\dots
\end{equation}
satisfies the {\it difference condition}
\begin{equation}\label{E:difference condition}
\sum_{m\in\mathcal Z}m\leq k
\end{equation}
for all downward paths $\mathcal Z$ in $\mathcal F\sp{(k_0, k_1, k_2)}$.
So, for example, $f_1$ in the first row must be $\leq k_2$ because of (\ref{E:difference condition}) for the downward path $\mathcal Z=\{f_1,0,k_1,0,k_0\}$.
We say that colored partitions (\ref{E:colored partition on the array}) with $(k_0, k_1, k_2)$-admissible arrays of frequencies (\ref{E:array of frequencies}) are {\it $(k_0, k_1, k_2)$-admissible colored partitions}.

As explained, we consider the natural numbers at different places in the array (\ref{E:C(1)2 principally specialized array}) as different. Likewise, in the array of frequencies (\ref{E:array of frequencies}) the entry $f_1$ in the first row denotes the frequency of the part $1$ in the first row, which may be different from the entry $f_1$ in the third row denoting the frequency of the part $1$ in the third row. However, sometimes we need to write explicitly the coloring of elements in the array (\ref{E:C(1)2 principally specialized array}) and change the notation of frequency arrays (\ref{E:array of frequencies}) and extended frequency arrays (\ref{E:extended array of frequencies}) accordingly, like in the following example:

\begin{example} \label{Ex:(2,0,0)-admissible colored partitions} Here we list all $(2,0,0)$-admissible colored partitions for $n\leq8$.
First we write explicitly one possible coloring of elements in the array (\ref{E:C(1)2 principally specialized array})
\begin{equation}\label{E:colored C(1)2 principally specialized array}
\begin{matrix}
1_{1}& & 3_{1}& & 5_{1}& & 7_{1}& &  9_{1}&  \\
& 2_{1}& &4_{1} & &6_{1} & & 8_{1} & & 10_{1} \\
1_{2}& &3_{2} & &5_{2} & &7_{2} & &9_{2} & \\
&2_{2} & &4_{2} & &6_{2} & &8_{2} & &10_{2}  \\
1_{3}& &3_{3} & &5_{3} & &7_{3} & & 9_{3} &
\end{matrix}\quad\dots \,.
\end{equation}
and change the notation for the extended array of frequencies (\ref{E:extended array of frequencies}) accordingly
\begin{equation}\label{E:colored extended array of frequencies}
\begin{matrix}
0&&f_{1_{1}}& & f_{3_{1}}& & f_{5_{1}}& & f_{7_{1}}&   \\
& 0&& f_{2_{1}}& &f_{4_{1}} & &f_{6_{1}} & & f_{8_{1}}  \\
0&&f_{1_{2}}& &f_{3_{2}} & &f_{5_{2}} & &f_{7_{2}} &  \\
& 0&&f_{2_{2}} & &f_{4_{2}} & &f_{6_{2}} & &f_{8_{2}}   \\
2&&f_{1_{3}}& &f_{3_{3}} & &f_{5_{3}} & &f_{7_{3} }&
\end{matrix}\quad\dots
\end{equation}
Then we have $(2,0,0)$-admissible colored partitions for $n\leq8$:
\begin{eqnarray}\label{colored patrition}
1 & = & 1_3\nonumber\\
2 & = & 2_2 = 1_3+1_3\nonumber\\
3 & = & 3_2 = 3_3 =  2_2 + 1_3\nonumber\\
4 & = & 4_1 = 4_2 = 3_2 + 1_3 = 3_3 + 1_3 = 2_2 + 2_2 \nonumber\\
5 & = & 5_1 = 5_2 = 5_3 = 4_1 + 1_3 = 4_2 + 1_3 = 3_2 + 2_2 = 3_3 + 2_2 = 3_3 + 1_3 + 1_3\nonumber\\
6 & = & 6_1 = 6_2 = 5_1 + 1_3 = 5_2 + 1_3 = 5_3 + 1_3 = 4_1 + 2_2 = 4_2 + 2_2 = 4_2 + 1_3 + 1_3\nonumber\\
7 & = & 7_1 = 7_2 = 7_3 = 6_1 + 1_3 = 6_2 + 1_3 = 5_1 + 2_2 = 5_2 + 2_2 = 5_3 + 2_2\nonumber\\
  & = & 5_2 + 1_3 + 1_3 = 5_3+ 1_3 + 1_3 = 4_1 + 3_2= 4_1 + 3_3 = 4_2 + 3_2 = 4_2 + 3_3  \nonumber\\
  & = & 4_2 + 2_2 +1_3 = 3_2 + 3_3 +1_3 = 3_3 + 3_3 + 1_3 \nonumber\\
8 & = & 8_1 = 8_2 = 7_1 + 1_3 = 7_2 + 1_3 = 7_3 +1_3 = 6_1 + 2_2 = 6_2 + 2_2 = 6_1 + 1_3 + 1_3 \nonumber\\
  & = & 6_2 + 1_3 + 1_3 =   5_1 + 3_2 = 5_1 + 3_3 = 5_2 + 3_2 = 5_2 + 3_3= 5_3 + 3_2 = 5_3 + 3_3 \nonumber\\
  & = & 5_2 + 2_2 + 1_3 = 5_3 + 2_2 + 1_3 = 4_1 + 4_1 = 4_1 + 4_2 = 4_2 + 4_2 = 4_1+3_3+1_3\nonumber\\
  & = & 4_2+3_3+1_3 = 4_2+3_2+1_3 = 4_2+2_2+2_2 = 3_3+3_3 +1_3+1_3\nonumber\ .
\end{eqnarray}
Note that for any $(2,0,0)$-admissible colored partition $n=\sum_{a\in\mathcal N}f_a\cdot a$ difference conditions (\ref{E:difference condition}) imply
$$
f_{1_{1}} = f_{2_{1}}= f_{3_{1}}=f_{1_{2}}=0.
$$
Also note that, for example, $8= 3_3+3_3 +1_3+1_3$ is a  $(2,0,0)$-admissible colored partition, and $8= 3_2+3_2 +1_3+1_3$ is not since  difference condition (\ref{E:difference condition}) is violated:
$$
f_{3_2}+f_{1_3}=2+2=4>2.
$$
\end{example}

\begin{remark} \label{R:(2,0,0)-admissible colored partitions}
The partitions in the above example are (essentially) the partitions in Example 3 in \cite{PS2}; only the coloring of the array (\ref{E:C(1)2 principally specialized array}) is different.
\end{remark}

We extend these notions for $C\sp{(1)}_\ell$, $\ell\geq1$, by starting with the array $\mathcal N=\mathcal N_{2\ell+1}=\mathcal N_{C^{(1)}_\ell}$ of $\ell$ copies of the set of natural numbers $\mathbb N$ and the additional set of odd numbers, arranged in $2\ell +1$ rows and diagonals of width $w=2\ell +1$; with numbers in $\mathcal N$ increasing by one going to the right on any diagonal. For example, for $\ell=3$ we have
the extended array of frequencies $\mathcal F\sp{(k_0, k_1, k_2, k_3)}$
\begin{equation}\label{E:extended array of frequencies w=7}
\begin{matrix}
 k_{3}&& f_1& &f_ 3& & f_5& & f_7&    \\
 & 0&&f_ 2& &f_4 & &f_6 & & f_{8}  \\
 k_{2}&& f_1& &f_ 3& & f_5& & f_7&    \\
 & 0&&f_ 2& &f_4 & &f_6 & & f_{8}  \\
 k_{1}&& f_1& &f_ 3& & f_5& & f_7&   \\
 & 0&&f_ 2& &f_4 & &f_6 & & f_{8}  \\
 k_{0}&& f_1& &f_ 3& & f_5& & f_7&    \\
\end{matrix}\quad\dots
\end{equation}
and the corresponding notion of $(k_0, k_1, k_2, k_3)$-admissible colored partitions on the array $\mathcal N_7$.

\begin{conjecture}\label{C:conjecture}  Let $\ell\geq2$. The principally specialized character
$$
\frac{\prod_{a\in\{0\}^\ell\cup D(k_0+1, k_1+1,\dots, k_\ell+1);\,
b\in \Delta( k_1+1,\dots, k_\ell+1);\,
j\equiv a,\pm b\mod (2\ell+2k+2)}(1-q^j)
}
{\prod_{j\, \text{odd}}(1-q^j)\prod_{j\in\mathbb N}(1-q^j)^\ell}
$$
is the generating function for the number of $(k_0, k_1, \dots, k_\ell)$-admissible partitions
$$
n=\sum_{a\in\mathcal N_{2\ell+1}}f_a\cdot a.
$$
\end{conjecture}

\begin{remark}\label{R:known w odd cases}
Conjecture \ref{C:conjecture} is true for $(1,0,\dots,0)$-admissible colored partitions and any $\ell\geq2$ (see \cite{PS1}).
In  \cite{PS2} Conjecture \ref{C:conjecture} is formulated for $(k,0,\dots,0)$-admissible colored partitions with $\ell\geq2$ and $k\geq2$, and it was checked for some small values of $\ell$, $k$ and $n$  (cf. \cite{P}). 

In \cite{MP} the product formula for the generating function  for the number of $(k_0, k_1)$-admissible partitions on the array $\mathcal N_3$ is given; in the $\ell=1$ case it is Lepowsky's product formula for $A_1^{(1)}$ and the key ingredient---a construction of combinatorial bases of standard $A_1^{(1)}$-modules---is independently obtained in \cite{FKLMM} and \cite{F}.

The combinatorial identities for $(1,0)$-admissible partitions and  for $(0,1)$-admissible partitions are equivalent\footnote{
In \cite{MP} certain spanning sets $B(\Lambda_0)\subset  L(\Lambda_0)$ and $B(\Lambda_1)\subset  L(\Lambda_1)$ of the two fundamental $A_1^{(1)}$-modules are constructed. The following three statements are equivalent: (i) $B(\Lambda_0)$ and $B(\Lambda_1)$ are linearly independent, (ii) the generating functions for $(1,0)$-admissible and $(0,1)$-admissible partitions are the principally specialized characters $\text{ch}_q L(\Lambda_0)$ and $\text{ch}_q L(\Lambda_1)$, and (iii) the two Capparelli identities hold.
 }
to the two Capparelli identities \cite{C}. Moreover, it seems that these identities are related to the purely combinatorial approach in \cite{AAG}, based on Capparelli's identity which was found using representation theory of $A_2^{(2)}$.
\end{remark}

\begin{example} \label{Ex:some colored partitions for w odd}

Conjecture \ref{C:conjecture} is based on a computer experiment.  Here we present some results for admissible partitions for $n$ up to $20$.

We write
$$
(1,0,1)\sim ( r\equiv 1,1,3,4,4,6,6,7,9,9 \mod10)
$$
if the number of $(1,0,1)$-admissible colored partitions of $n$  is equal to the number of colored partitions of $n$ with parts $ r\equiv 1,1,3,4,4,6,6,7,9,9 \mod10$ for $n\leq 20$  (here $ r\equiv 1,1 \mod10$ means that parts $ r\equiv 1 \mod10$ come in two colors). Since the parameters $(k_0,k_1,k_2)$ and $(k_2,k_1,k_0)$ give ``isomorphic'' colored partitions, we list below only mutually different conjectured identities.
\smallskip

\noindent
For $w=5$ we have:
$$
{\allowdisplaybreaks
\begin{aligned}
(1,0,0)&\sim (r \ \text{odd}; \    r\equiv 4 \mod 8),\\
(2,0,0)&\sim (r \ \text{odd}; \    r\equiv 2,4,5,6,8 \mod 10),\\
(1,1,0)&\sim ( r \ \text{odd}; \   r\equiv1,3,5,7 ,9 \mod 10),\\
(1,0,1)&\sim ( r\equiv 1,1,3,4,4,6,6,7,9,9 \mod10),\\
(0,2,0)&\sim ( r\equiv1,2,2,3,3,7,7,8,8,9  \mod 10),\\
(3,0,0)&\sim ( r \ \text{odd}; \   r\equiv 2,3,4,5,6,7,8,9,10  \mod 12),\\
(2,1,0)&\sim ( r \ \text{odd}; \  r\equiv 1,2,4,5,6,7,8,10,11 \mod 12),\\
(2,0,1)&\sim ( r \ \text{odd}; \  r\equiv 1,2,4,5,6,7,8,10,11 \mod12),\\
(1,2,0)&\sim ( r \ \text{odd}; \  r\equiv 1,2,3,4,6,8,9,10,11  \mod 12),\\
(0,3,0)&\sim  ( r \ \text{odd}; \ r\equiv 2,3,4,5,6,7,8,9,10 \mod 12),\\
\end{aligned}
}
$$
$$
{\allowdisplaybreaks
\begin{aligned}
(3,0,1)&\sim (r \ \text{odd}; \   r\equiv 1,2,3,4,6,6,7,8,8,10,11,12,13\mod14),\\
(2,1,1)&\sim (r, r \ \text{odd}; \  r\equiv 1,4,6,8,10,13  \mod 14),\\
(2,0,2)&\sim (r\equiv 1,1,2,2,3,5,5,5,6,6,8,8,9,9,9,11,12,12,13,13 \mod14),\\
(0,4,0)&\sim (r \equiv 1,2,2,3,3,3,4,4,5,5,9,9,10,10,11,11,11,12,12,13\mod14),\\
(3,0,2)&\sim (r, r \ \text{odd}; \   r\equiv 2,2,6,6,8,10,10,14,14\mod16).\\
\end{aligned}
}
$$
For $w=7$ we have:
$$
\begin{aligned}
(1,0,0,1)&\sim (r \ \text{odd}; \ r\equiv 1,3,4,5,7,8,9,11 \mod12),\\
(0,1,1,0)&\sim (r \ \text{odd}; \ r\equiv 1,2,3,5,7,9,10,11 \mod12),\\
(2,0,0,1)&\sim (r \ \text{odd}; \ r\equiv 1,2,3,4,5,6,6,8,8,9,10,11,12,13 \mod14),\\
(2,0,0,2)&\sim (r, r \ \text{odd}; \ r\equiv 2,2,4,6,6,6,10,10,10,12,14,14 \mod16).
\end{aligned}
$$
For $w=9$ we have:
$$
\begin{aligned}
(1,0,0,0,1)\sim &\ (r\equiv 1,1,3,3,4,4,5,6,6,8,8,9,10,10,11,11,13,13 \mod14),\\
(0,1,1,0,0)\sim &\ (r \ \text{odd}; \ r\equiv 1,2,3,4,5,7,9,10,11,12,13 \mod14),\\
(0,1,1,0,1)\sim &\ (r \ \text{odd}; \ r\equiv 1,1,2,3,4,4,6,6,7,8,\\
&9,10,10,12,12,13,14,15\mod16),\\
(2,1,0,0,1)\sim &\ (r, r \ \text{odd}; \ r\equiv 1,2,4,4,5,6,7,8,8,\\
&10,10,11,12,13,14,14,16,17\mod18).\\
\end{aligned}
$$
\end{example}

\begin{example} \label{Ex 3.6:some colored partitions for w odd}
In the following identities, the right hand sides can be interpreted as the number of colored partitions satisfying congruence conditions, with the extra requirement that in one color the parts have to be distinct.

Associated to the module $L(0,1,0)$ we have that, for $n \le 20$, the number of $(0,1,0)$-admissible colored partitions of $n$ equals
  \begin{equation*}
    \operatorname{coeff}_{q^n}
\prod_{r \equiv 2 \mod 4} (1+q^r) \big / \prod_{r \equiv 1,3,5,7 \mod 8}(1-q^r).
  \end{equation*}
  
Associated to the module $L(1,1,1)$ we have that, for $n \le 20$,  the number of $(1,1,1)$-admissible colored partitions of $n$ equals
\begin{equation*}
\operatorname{coeff}_{q^n}
\prod_{r \text{odd}} (1+q^r) \big / \prod_{r \text{odd}}(1-q^r)^2.
\end{equation*}

\end{example}

\section{Arrays with even width  $w\geq4$} \label{Even array width w}

Let $\mathcal N^e=\mathcal N_4$ be the array of natural numbers
\begin{equation}\label{E:w even array}
\begin{matrix}
1& & 3& & 5& & 7& &  {9}&  \\
& 2& &4 & &6 & & {8} & & 10 \\
1& &3 & &5 & &7 & &{9} & \\
&2 & &4 & &6 & &{8} & &10  \\
\end{matrix}\quad\dots \,.
\end{equation}
This array consists of two copies of the set of natural numbers $\mathbb N$ arranged in $4$ rows, with diagonals of width $w=4$. We consider colored partitions
\begin{equation}\label{E:colored partition on w even}
n=\sum_{a\in\mathcal N^e}f_a\cdot a,
\end{equation}
where $f_a$ is the frequency of the part $a\in\mathcal N^e$ in the colored partition (\ref{E:colored partition on w even}) of $n$.
Let $k_0, k_1, k_2\in\mathbb N_0$, $k=k_0+k_1+k_2>0$. We say that an array of frequencies    $\mathcal F^e$
\begin{equation}\label{E:array of frequencies on w even}
\qquad\qquad
\begin{matrix}
f_1& &f_ 3& & f_5& & f_7&    \\
&f_ 2& &f_4 & &f_6 & & f_{8}  \\
f_1& &f_ 3& & f_5& & f_7&    \\
&f_ 2& &f_4 & &f_6 & & f_{8}  \\
\end{matrix}\quad\dots
\end{equation}
{\it is $(k_0, k_1, k_2)^e$-admissible} if {\it the extended array of frequencies} $\mathcal F\sp{e(k_0, k_1, k_2)}$
\begin{equation}\label{E:the extended array of frequencies on w even}
\begin{matrix}
 k_{2}&& f_1& &f_ 3& & f_5& & f_7&    \\
 & 0&&f_ 2& &f_4 & &f_6 & & f_{8}  \\
 k_{1}&& f_1& &f_ 3& & f_5& & f_7&   \\
 &  k_{0}&&f_ 2& &f_4 & &f_6 & & f_{8}  \\
\end{matrix}\quad\dots
\end{equation}
satisfies the {\it difference condition}
\begin{equation}\label{E:difference condition on w even}
\sum_{m\in\mathcal Z}m\leq k
\end{equation}
for all downward paths $\mathcal Z$ in $\mathcal F\sp{e(k_0, k_1, k_2)}$.
So, for example, $f_1$ in the first row must be $\leq k_2$ because of (\ref{E:difference condition on w even}) for downward path $\mathcal Z=\{f_1,0,k_1,k_0\}$.
We say that colored partitions (\ref{E:colored partition on w even}) with $(k_0, k_1, k_2)^e$-admissible arrays of frequencies (\ref{E:array of frequencies on w even}) are {\it $(k_0, k_1, k_2)^e$-admissible colored partitions}.

We extend these notions for $\ell\geq1$ by starting with the array $\mathcal N^e=\mathcal N_{2\ell}$ of $\ell$ copies of the set of natural numbers $\mathbb N$  arranged in $2\ell$ rows and diagonals of width $w=2\ell$; with numbers in $\mathcal N^e$ increasing by one going to the right on any diagonal. For example, for $\ell=3$ ($w=6$) instead of (\ref{E:the extended array of frequencies on w even}) we have the extended array of frequencies
$\mathcal F\sp{e(k_0, k_1, k_2,k_3)}$
\begin{equation}\label{E:the extended array of frequencies on w even l=3}
\begin{matrix}
 k_{3}&& f_1& &f_ 3& & f_5& & f_7&    \\
 & 0&&f_ 2& &f_4 & &f_6 & & f_{8}  \\
 k_{2}&& f_1& &f_ 3& & f_5& & f_7&    \\
 & 0&&f_ 2& &f_4 & &f_6 & & f_{8}  \\
 k_{1}&& f_1& &f_ 3& & f_5& & f_7&   \\
 &  k_{0}&&f_ 2& &f_4 & &f_6 & & f_{8}  \\
\end{matrix}\quad\dots
\end{equation}
and the corresponding notion of $(k_0, k_1, k_2, k_3)^e$-admissible colored partitions on the array $\mathcal N_6$.

\begin{conjecture}\label{C:conjecture w even}  Let $k_0, k_1, \dots, k_\ell\in\mathbb N_0$, $k=k_0+k_1+\dots+k_\ell>0$. Then the generating function for the number of $(k_0, k_1, \dots, k_\ell)^e$-admissible colored partitions
$$
n=\sum_{a\in\mathcal N_{2\ell}}f_a\cdot a
$$
is the infinite periodic product
\begin{equation}\label{E:product formula for w even}
\frac{\prod_{a\in\{0\}^\ell;\,
b\in \Delta( k_1+1,\dots, k_\ell+1);\,
j\equiv a,\pm b\mod (2\ell+2k+1)}
(1-q^j)}
{\prod_{j\in\mathbb N}(1-q^j)^\ell}.
\end{equation}
\end{conjecture}

\begin{example}
By using Example \ref{Ex:congruence triangle (3,2,1,1,2)} for the congruence triangle $\Delta(2,1,1,2)$, we see that the conjectured product (\ref{E:product formula for w even}) for $(2,1,0,0,1)^e$-admissible colored partitions is
$$
\prod_{j\equiv  1,1,2,3,3,4,4,5,5,5,6,6,7,7,8,8,9,9,10,10,11,11,12,12,12,13,13,14,14,15,16,16
\mod 17}(1-q^j)^{-1}.
$$
\end{example}

\begin{remark}\label{R: Rogers-Ramnujan identities}
Conjecture \ref{C:conjecture w even} is true for $\ell=1$: the extended array of frequencies
$\mathcal F\sp{e(k_0, k_1)}$ is
\begin{equation}\label{E:the extended array of frequencies on w even l=1}
\begin{matrix}
 k_{1}&& f_1& &f_ 3& & f_5& & f_7&   \\
 &  k_{0}&&f_ 2& &f_4 & &f_6 & & f_{8}  \\
\end{matrix}\quad\dots \,,
\end{equation}
and $(k_0,k_1)^e$-admissible  colored partitions are classical partitions
$$
n=\sum_{a\in\mathbb N}f_a\cdot a
$$
satisfying difference and initial conditions
$$
f_a+f_{a+1}\leq k,\qquad f_1\leq k_1.
$$
So the generating functions for the number of  $(1,0)^e$-admissible and $(0,1)^e$-admissible partitions
are the product sides of two Rogers-Ramanujan identities, and for $k=k_0+k_1>1$ we have the product sides of Gordon identities (cf. \cite{A1}, \cite{A2}, \cite{G}). These combinatorial identities have a Lie-theoretic interpretation (see \cite{LW}), but for $2\ell>2$ there is no obvious connection of
$(k_0, k_1, \dots, k_\ell)^e$-admissible colored partitions (\ref{E:colored partition on w even}) with  representation theory of affine Lie algebras.
\end{remark}

\begin{example} \label{Ex:some colored partitions for w even}
In this example we write
$$
(0,0,1)^e\sim ( r\equiv 1,3,4,6 \mod 7)
$$
if the number of $(0,0,1)^e$-admissible colored partitions of $n$  is equal to the number of colored partitions of $n$ with parts $ r\equiv 1,3,4,6 \mod 7$ for $n\leq 20$.

For $w=2$ we have Rogers-Ramanujan and Gordon identities:
$$
\begin{aligned}
(1,0)^e&\sim (r\equiv 2,3 \mod{5}),\\
(0,1)^e&\sim (r\equiv 1,4 \mod{5}),\\
(2,0)^e&\sim (r\equiv 2,3,4,5 \mod{7}),\\
(1,1)^e&\sim (r\equiv 1,3,4,6 \mod{7}),\\
(0,2)^e&\sim (r\equiv 1,2,5,6 \mod{7}),\\
(3,0)^e&\sim (r\equiv 2,3,4,5,6,7 \mod{9}),\\
(2,1)^e&\sim (r\equiv 1,3,4,5,6,8 \mod{9}),\\
(1,2)^e&\sim (r\equiv 1,2,4,5,7,8 \mod{9}),\\
(0,3)^e&\sim (r\equiv 1,2,3,6,7,8 \mod{9}).\\
\end{aligned}
$$

For $w=4$ we have:
$$
\begin{aligned}
(1,0,0)^e&\sim (r\equiv 2,3,4,5 \mod{7}),\\
(0,1,0)^e&\sim (r\equiv 1,2,5,6 \mod{7}),\\
(0,0,1)^e&\sim (r\equiv 1,3,4,6 \mod{7}),\\
(2,0,0)^e&\sim (r\equiv  2,3,4,4,5,5,6,7 \mod{9}),\\
(1,1,0)^e&\sim (r\equiv 1,2,3,4,5,6,7,8  \mod{9}),\\
(1,0,1)^e&\sim (r\equiv 1,2,3,4,5,6,7,8 \mod{9}),\\
(0,2,0)^e&\sim (r\equiv 1,2,2,3,6,7,7,8 \mod{9}),\\
(0,1,1)^e&\sim (r\equiv 1,1,3,4,5,6,8,8 \mod{9}),\\
(0,0,2)^e&\sim (r\equiv 1,2,3,4,5,6,7,8 \mod{9}),\\
\end{aligned}
$$
$$
\begin{aligned}
(3,0,0)^e&\sim (r\equiv 2,3,4,4,5,5,6,6,7,7,8,9 \mod{11}),\\
(2,1,0)^e&\sim (r\equiv 1,2,3,4,5,5,6,6,7,8,9,10 \mod{11}),\\
(2,0,1)^e&\sim (r\equiv  1,2,3,4,4,5,6,7,7,8,9,10\mod{11}),\\
(1,2,0)^e&\sim (r\equiv 1,2,2,3,4,5,6,7,8,9,9,10 \mod{11}),\\
(1,1,1)^e&\sim (r\equiv 1,1,3,3,4,5,6,7,8,8,10,10 \mod{11}),\\
(1,0,2)^e&\sim (r\equiv  1,2,2,3,5,5,6,6,8,9,9,10\mod{11}),\\
(0,3,0)^e&\sim (r\equiv 1,2,2,3,3,4,7,8,8,9,9,10 \mod{11}),\\
(0,2,1)^e&\sim (r\equiv 1,1,2,3,4,5,6,7,8,9,10,10 \mod{11}),\\
(0,1,2)^e&\sim (r\equiv 1,1,2,4,4,5,6,7,7,9,10,10 \mod{11}),\\
(0,0,3)^e&\sim (r\equiv 1,2,3,3,4,5,6,7,8,8,9,10 \mod{11}).\\
\end{aligned}
$$

For $w=6$ we have:
$$
\begin{aligned}
(1,0,0,0)^e&\sim  (r\equiv 2,3,4,5,6,7 \mod{9}),\\
(0,1,0,0)^e&\sim  (r\equiv 1,2,4,5,7,8 \mod{9}),\\
(0,0,1,0)^e&\sim  (r\equiv 1,2,3,6,7,8 \mod{9}),\\
(0,0,0,1)^e&\sim  (r\equiv 1,3,4,5,6,8 \mod{9}),\\
(2,0,0,0)^e&\sim  (r\equiv 2,3,4,4,5,5,6,6,7,7,8,9 \mod{11}),\\
(1,1,0,0)^e&\sim  (r\equiv 1,2,3,4,4,5,6,7,7,8,9,10 \mod{11}).\\
\end{aligned}
$$

For $w=8$ we have:

$(0,1,1,0,1)^e\sim $
$(r\equiv 1,1,1,2,3,4,4,4,6,6,6,7,\\8,9,9,9, 11,11,11,12,13,14,14,14 \mod{15})$,
\smallskip

$(2,1,0,0,1)^e\sim$
$(r\equiv 1,1,2,3,3,4,4,5,5,5,6,6,7,7,8,8,\\9,9,10,10,11,11,12,12,12,13,13,14,14,15,16,16 \mod 17)$,
\smallskip

$(0,1,1,1,1)^e\sim $
$(r\equiv 1,1,1,1,3,3,3,4,5,5,5,6,7,7,8,8,\\9,9, 10,10,11,12,12,12,13,14,14,14,16,16,16,16\mod{17})$.
\end{example}

\begin{remark}
It seems that for pairs (level $k$, rank $\ell$)  there is some sort of duality
$$
(k, \ell)\longleftrightarrow (\ell, k).
$$
In particular, the Rogers-Ramanujan case $k=1$, $w=2$ is  self-dual and $k=2$, $w=2$ is dual to $k=1$, $w=4$. In the self-dual case $k=2$, $w=4$ we see that $(1,1,0)^e$-admissible,  $(1,0,1)^e$-admissible and  $(0,0,2)^e$-admissible partitions have the same product formula, but already for $n=1$ and $n=2$ we see that three types of colored partitions on $\mathcal N_4$ are mutually different.
\end{remark}

\section{An algorithm for constructing admissible arrays of frequencies}\label{An algorithm}

In this section we describe an algorithm for constructing $(k_0, k_1,\dots, k_\ell)$-admissible and $(k_0, k_1,\dots, k_\ell)^e$-admissible arrays of frequencies.
First we consider the simplest case of $(k, 0, 0)$-admissible arrays of frequencies    $\mathcal F$. The difference condition (\ref{E:difference condition}) forces
the extended array of frequencies $\mathcal F\sp{(k,0,0)}$ to look like
\begin{equation}\label{E:(k,0,0)-extended array of frequencies}
\begin{matrix}
 0&& 0& &0& & f_{5_1}& & f_{7_1}&    \\
 & 0&&0& &f_{4_1} & &f_{6_1} & & f_{8_1}  \\
 0&& 0& &f_ {3_2}& & f_{5_2}& & f_{7_2}&   \\
 & 0&&f_ {2_2}& &f_{4_2} & &f_{6_2} & & f_{8_2}  \\
 k&& f_{1_3}& &f_ {3_3}& & f_{5_3}& & f_{7_3}&    \\
\end{matrix}\dots\,.
\end{equation}
In order to avoid any confusion and facilitate the exposition, here we have colored our arrays as in Example  \ref{Ex:(2,0,0)-admissible colored partitions}.

At the beginning we let $\mathcal F=0$, i.e. all the frequencies in $\mathcal F$ are zero. Then we construct a nontrivial $(k, 0, 0)$-admissible array of frequencies by changing $\mathcal F$ in steps. We start our construction from the top $ {5_1}$ till the bottom ${1_3}$ of the first diagonal $ \{{5_1},  {4_1},  {3_2},  {2_2},  {1_3}\}\subset\mathcal N $ of lenght $5$, then from the top till the bottom of the second diagonal $  \{{7_1},  {6_1},  {5_2},  {4_2},  {3_3}\}$ of lenght $5$, and so on. At each step for $a\in\mathcal N$ we choose a frequency $f_a$ and determine the corresponding maximum
\begin{equation}\label{E:maxima at a}
m_a=\max\,\big\{\sum_{c\in\mathcal Z} f_c\mid \mathcal Z\text{ \ is a downward path which ends in \ } a\big\}.
\end{equation}

So we start with $a=5_1$ and we choose any value $f_{5_1}\in\{0, 1,\dots , k\}$. For  the first point $5_1$ there is only one downward path $\mathcal Z=\{5_1\}$ which ends in $5_1$, so $$m_{5_1}=f_{5_1}.$$ 
For the second point $4_1$ we should choose $f_{4_1}$. There are two downward paths which end in $4_1$:  $\mathcal Z_1=\{5_1, 4_1\}$ and  $\mathcal Z_2=\{3_1, 4_1\}$. Since $f_{3_1}=0$, we have
$$m_{4_1}=f_{5_1}+f_{4_1}$$ 
and the level $k$ difference condition $m_{4_1}\leq k$ forces us to choose 
$$f_{4_1}\leq k-f_{5_1}.$$
At each step we choose $f_a$ and determine $m_a$, so assume we completed a list of frequencies $f_a$ and maxima $m_a$ in the first two diagonals and in the top two places ${9_1},  {8_1}$ in the third diagonal. Next we should choose $f_{b}$ for $b=7_2$. Denote the points $ {6_1}$ and $ {8_1}$---the points above and adjacent to $b$---as $x$ and $y$:
 \begin{equation}\label{E:construction of (k,0,0)-admissible array of frequencies}
\begin{matrix}
& & & & \bullet& & \bullet& &  \bullet&   \\
& & &\bullet& &x & & y & &  \\
& &\bullet & &\bullet & &b& & & \\
&\bullet & &\bullet & & & & & &  \\
\bullet& &\bullet & & & & & &  &
\end{matrix}
\end{equation}
Then choose a frequency $f_b$ so that $f_b+max\{m_x, m_y\}\leq k$. Since every downward path $\mathcal Z$ comes to $b$ via $x$ or via $y$, it is clear that
\begin{equation}\label{E:recursive condition for frequencies}
m_b=f_b+max\{m_x, m_y\}\leq k.
\end{equation}

In this sequence of steps we have constructed the frequencies $f_c$ for $c$ denoted as $x, y, b$ or $\bullet$ in (\ref{E:construction of (k,0,0)-admissible array of frequencies}), and all the other frequencies in $\mathcal F$ are (still) $0$. By construction, i.e. by the recursive condition (\ref{E:recursive condition for frequencies}) for each constructed $f_c$, the difference condition (\ref{E:difference condition}) holds for all downward paths $\mathcal Z$ in $\mathcal F\sp{(k,0,0)}$. Hence the (so far) constructed sequence of frequencies is $(k,0,0)$-admissible. 
In this way we proceed till the end of a finite (chosen in advance) number of diagonals.

\begin{remark}
Assume we want to construct and count all colored partitions
$$
n=\sum_{a\in\mathcal N_5}f_a\cdot a\quad\text{for}\quad n\leq 15,
$$
with $(k, 0, 0)$-admissible frequencies $f_a$. Since the biggest part in a partition of $n\leq 15$ is at most $15$, it is enough to consider only the first eight nontrivial diagonals of frequency array
(\ref{E:(k,0,0)-extended array of frequencies}).
We can write compactly
the first eight full-length diagonals in the array (\ref{E:C(1)2 principally specialized array}), the constructed part of the frequency array
(\ref{E:construction of (k,0,0)-admissible array of frequencies}), and the constructed part of the maxima array as\footnote{
The diagram (\ref{E:construction for  (k,0,0)-admissible compactly}) is obtained from (\ref{E:(k,0,0)-extended array of frequencies}) by a (clockwise) 45 degree rotation, followed by the shear linear transformation $(x,y) \mapsto (x+y,y)$. In particular,  the rows in (\ref{E:construction for  (k,0,0)-admissible compactly}) are obtained from diagonals in (\ref{E:(k,0,0)-extended array of frequencies}).
}
\begin{equation}\label{E:construction for  (k,0,0)-admissible compactly}
\begin{matrix}
1_3&2_2 &3_2 &4_1 & 5_1 \\
3_3&4_2&5_2&6_1&7_1\\
5_3&6_2&7_2&8_1&9_1\\
7_3&8_2&9_2&10_1&11_1\\
9_3&10_2&11_2&12_1&13_1\\
11_3&12_2&13_2&14_1&15_1\\
13_3&14_2&15_2&16_1&17_1\\
15_3&16_2&17_2&18_1&19_1\\
\end{matrix}\quad
\begin{matrix}
f_{1_3}&f_{2_2}&f_{3_2} &f_{4_1} &f_ {5_1} \\
f_{3_3}&f_{4_2}&f_{5_2}&f_{6_1}&f_{7_1}\\
\cdot&\cdot&\cdot&f_{8_1}&f_{9_1}\\
\cdot&\cdot&\cdot&\cdot&\cdot\\
\cdot&\cdot&\cdot&\cdot&\cdot\\
\cdot&\cdot&\cdot&\cdot&\cdot\\
\cdot&\cdot&\cdot&\cdot&\cdot\\
\cdot&\cdot&\cdot&\cdot&\cdot\\
\end{matrix}
\quad
\begin{matrix}
m_{1_3}&m_{2_2}&m_{3_2} &m_{4_1} &m_ {5_1} \\
m_{3_3}&m_{4_2}&m_{5_2}&m_{6_1}&m_{7_1}\\
\cdot&\cdot&\cdot&m_{8_1}&m_{9_1}\\
\cdot&\cdot&\cdot&\cdot&\cdot\\
\cdot&\cdot&\cdot&\cdot&\cdot\\
\cdot&\cdot&\cdot&\cdot&\cdot\\
\cdot&\cdot&\cdot&\cdot&\cdot\\
\cdot&\cdot&\cdot&\cdot&\cdot\\
\end{matrix}
\end{equation}
Then we can choose a frequency $0\leq f_{7_2}\leq k$ for $7_2$  so that the recursive condition (\ref{E:recursive condition for frequencies}) holds, i.e.
\begin{equation}\label{E:recursive relation for frequencies example}
m_{7_2}=f_{7_2}+max\{m_{6_1}, m_{8_1}\}\leq k.
\end{equation}
In this way we proceed till the end of the eighth row.
\end{remark}
Now we consider, again for $\ell=2$, a construction of
$(k_0,k_1,k_2)$-admissible or $(k_0,k_1,k_2)^e$-admissible arrays of frequencies $\mathcal F$ or $\mathcal F^e$ with non-zero frequencies in a finite (chosen in advance) number of diagonals.
We denote with $\circ$ the places of the array of not yet constructed frequencies in the chosen finite  number of diagonals as
$$
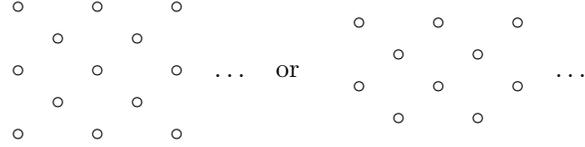
\begin{matrix}
\circ&&\circ&&\circ&\\
&\circ&&\circ&&\\
\circ&&\circ&&\circ& \\
 &\circ&&\circ&&  \\
\circ &&\circ&&\circ&
\end{matrix}
\dots\quad\text{or}\qquad
\begin{matrix}
\circ&&\circ&&\circ&\\
&\circ&&\circ&&\\
\circ&&\circ&&\circ& \\
 &\circ&&\circ&&  \\
\end{matrix}\dots
$$
and we extend them on the left with the prescribed fixed frequencies
$$
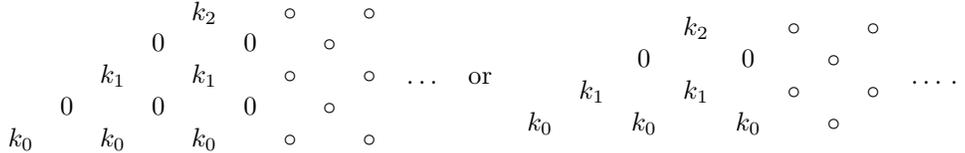
\begin{matrix}
&&&&k_2&&\circ&&\circ&\\
&&&0&&0&&\circ&&\\
& & k_1& & k_1& &\circ& &\circ & \\
&0 & &0& &0 & &\circ & &  \\
k_0& &k_0 & &k_0 & &\circ & & \circ &
\end{matrix}\dots\quad\text{or}\
\begin{matrix}
&&&&k_2&&\circ&&\circ&\\
&&&0&&0&&\circ&&\\
& & k_1& & k_1& &\circ& &\circ & \\
&k_0 & &k_0& &k_0 & &\circ & &  \\
\end{matrix}\dots \,.
$$
Note that added fixed frequencies satisfy the difference conditions (\ref{E:difference condition}) or
(\ref{E:difference condition on w even}).
Now the first diagonal is fixed and ``constructed'',  and the corresponding diagonal of the maxima is
$$
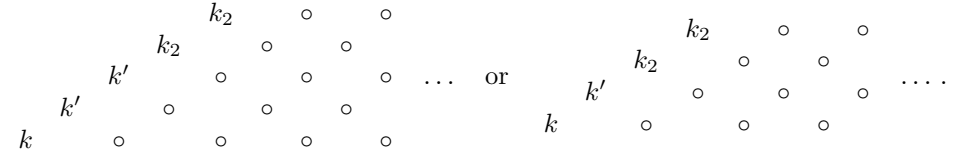
\begin{matrix}
&&&&k_2&&\circ&&\circ&\\
&&&k_2&&\circ&&\circ&&\\
& & k'& & \circ& &\circ& &\circ & \\
&k'& &\circ& &\circ & &\circ & &  \\
k& &\circ& &\circ & &\circ & & \circ &
\end{matrix}\dots\quad\text{or}\
\begin{matrix}
&&&&k_2&&\circ&&\circ&\\
&&&k_2&&\circ&&\circ&&\\
& & k'& & \circ& &\circ& &\circ & \\
&k& &\circ& &\circ & &\circ & &  \\
\end{matrix}\dots\,.
$$
with $k'=k_1+k_2$ and $k=k_0+k_1+k_2$. So we pass to the first point in the second diagonal, choose a frequency from the set $\{0, 1,\dots , k\}$ and determine the corresponding maximum. After that we pass to the second point and so on, just like before, except that in the second diagonal four or three frequencies are already fixed and if the condition (\ref{E:difference condition}) or (\ref{E:difference condition on w even}) does not hold, we should return to the beginning of the diagonal and give another try with the frequency of the first point. In the third diagonal we have to construct frequency for three points---the other two or one are already given and fixed, and so on. It is clear that the constructed sequence of frequencies is $(k_0,k_1,k_2)$-admissible.

In the Python code below we consider the $w$ odd and $w$ even cases simultaneously with the general initial conditions for parameters $[k_1, k_2, \dots, k_{w-1}, k_w]$: for a frequency array $\mathcal F$ for $\mathcal N_w$ we have the extended frequency array
$\mathcal F^{[k_1, k_2, \dots, k_{w-1}, k_w]}$
$$
\begin{matrix}
&&&&k_w&&\circ&&\circ&\\
&&& k_{w-1}&& k_{w-1}&&\circ&&\\
& &  k_{w-2}& &  k_{w-2}& &\circ& &\circ & \dots\,.\\
& k_{w-3} & & k_{w-3}& & k_{w-3}& &\circ & &  \\
& && &\vdots & & &\vdots &  &
\end{matrix}
$$
and the corresponding notion of  $[k_1, k_2, \dots, k_{w-1}, k_w]$-admissible colored partitions on the array $\mathcal N_w$. Note that
\begin{itemize}
\item  $(k_0,k_1,\dots, k_\ell)$-admissible $=$ $[k_0,0,k_1,0,\dots,0, k_\ell]$-admissible for $w=2\ell+1$,
\item  $(k_0,k_1,\dots, k_\ell)^e$-admissible $=$ $[k_0,k_1,0,\dots,0, k_\ell]$-admissible for $w=2\ell$.
\end{itemize}
We conjecture product formulas for generating  functions for the number of  the corresponding admissible colored partitions, and it seems that for any other  $[k_1, \dots, k_w]$
there is no infinite periodic product formula for the
generating function for the number of the corresponding $[k_1, \dots, k_w]$-admissible colored partitions.

\begin{remark}\label{R:IC for C(1)4 compactly} In higher ranks it may be convenient to write the extended $[k_1, \dots, k_w]$-admissible  frequency arrays compactly as infinite matrices
\begin{equation}\label{E:infinite matrix of frequencies}
F=(f_{ij})_{i=0,1,\dots; \ j=w, \dots, 1}
\end{equation}
with fixed prescribed frequencies in the upper left corner of the matrix, depending on
odd $w=2\ell+1$ or even $w=2\ell$. Here $j=w$ denotes the left column, and $j=1$ the right column.
For  $\ell=4$ we have matrices of the form
$$
\begin{matrix}
k_1&k_2&k_3&k_4&k_5&k_6&k_7&k_8&k_9\\
k_1&k_2&k_3&k_4&k_5&k_6&k_7&k_8&\bullet\\
k_1&k_2&k_3&k_4&k_5&k_6&\bullet&\bullet&\bullet\\
k_1&k_2&k_3&k_4&\bullet&\bullet&\bullet&\bullet&\bullet\\
k_1&k_2&\bullet&\bullet&\bullet&\bullet&\bullet&\bullet&\bullet\\
\circ&\circ&\circ&\circ&\circ&\circ&\circ&\circ&\circ\\
\circ&\circ&\circ&\circ&\circ&\circ&\circ&\circ&\circ\\
\circ&\circ&\circ&\circ&\circ&\circ&\circ&\circ&\circ\\
&&&&\vdots&&&&
\end{matrix}\qquad\text{or}\qquad
\begin{matrix}
k_1&k_2&k_3&k_4&k_5&k_6&k_7&k_8\\
k_1&k_2&k_3&k_4&k_5&k_6&k_7&\bullet\\
k_1&k_2&k_3&k_4&k_5&\bullet&\bullet&\bullet\\
k_1&k_2&k_3&\bullet&\bullet&\bullet&\bullet&\bullet\\
k_1&\bullet&\bullet&\bullet&\bullet&\bullet&\bullet&\bullet\\
\circ&\circ&\circ&\circ&\circ&\circ&\circ&\circ\\
\circ&\circ&\circ&\circ&\circ&\circ&\circ&\circ\\
\circ&\circ&\circ&\circ&\circ&\circ&\circ&\circ\\
&&&\vdots&&&&
\end{matrix}\ .
$$
For a frequency matrix (\ref{E:infinite matrix of frequencies}) we assume that there are finitely many non-zero elements $f_{ij}$ and we say that $F$ has a finite support. We can write the matrix $F$ as an
infinite sequence of rows
\begin{equation}\label{E:infinite matrix of frequencies2}
F=(\phi_{i})_{i=0,1,\dots},\qquad \phi_{i}=(f_{ij})_{ j=w, \dots, 1}.
\end{equation}
Like in (\ref{E:construction for  (k,0,0)-admissible compactly}), for a frequency matrix $F$ we have the associated {\it maxima matrix}
\begin{equation}\label{E:maxima matrix of F}
M=(m_{ij})_{i=0,1,\dots; \ j=w, \dots, 1},
\end{equation}
defined for $i=0$ as
$$
m_{0j}=k_w+\dots+k_{w-j+1},\quad j=1,\dots, w,
$$
and for $  i=1,2, \dots$ recursively like in (\ref{E:recursive relation for frequencies example}),
\begin{equation}\label{E:maxima matrix of F}
m_{i1}=f_{i1}, \quad
m_{ij}=f_{ij}+\max \{m_{i-1,j-1}, m_{i,j-1}\}, \quad    j=2, \dots,w.
\end{equation}
Our argument above shows that a frequency matrix $F$ is $[k_1, \dots, k_w]$-admissible if and only if
\begin{equation}\label{E:maxima matrix criteria}
m_{ij}\leq k\quad\text{for all} \quad  i=1,2,\dots, \quad j=1, \dots, w.
\end{equation}

We can also write our arrays of natural numbers $\mathcal N_2$, $\mathcal N_3$, $\mathcal N_4$, $\mathcal N_5$, \dots as infinite matrices $N_w^0$, $w=2, 3, 4, 5, \dots$ of natural numbers, extended with fixed zeros
in the upper left corner of the matrix,
$$
N_w^0=(n_{ij})_{i=0,1,\dots; \ j=w, \dots, 1},\qquad n_{ij}=\max\{0,2i-j\},
$$
i.e. as infinite matrices $N_2^0$, $N_3^0$, $N_4^0$, $N_5^0$, \dots
$$
\begin{matrix}
0&0\\
0&1\\
2&3\\
4&5\\
6&7\\
8&9\\
10&11\\
\vdots&&
\end{matrix},\qquad
\begin{matrix}
0&0&0\\
0&0&1\\
1&2&3\\
3&4&5\\
5&6&7\\
7&8&9\\
9&10&11\\
&\vdots&&
\end{matrix},\qquad
\begin{matrix}
0&0&0&0\\
0&0&0&1\\
0&1&2&3\\
2&3&4&5\\
4&5&6&7\\
6&7&8&9\\
8&9&10&11\\
&&\vdots&&
\end{matrix},\qquad
\begin{matrix}
0&0&0&0&0\\
0&0&0&0&1\\
0&0&1&2&3\\
1&2&3&4&5\\
3&4&5&6&7\\
5&6&7&8&9\\
7&8&9&10&11\\
&&\vdots&&
\end{matrix},\quad\dots\,.
$$
Then for each $[ k_1, \dots, k_w]$-admissible frequency matrix $F$ we have a $[ k_1, \dots, k_w]$-admissible  colored partition
$$
n=\sum_{i=0}^\infty\sum_{j=1}^w f_{ij}n_{ij}.
$$
\end{remark}

\begin{example}\label{Ex: a construction of (0,1,0)-admissible partitions} Here we describe a construction of all  $(0,1,0)$-admissible frequencies with the support in the first $5$ diagonals in the array $\mathcal N_5$. By using the notation in Remark \ref{R:IC for C(1)4 compactly}, we want to construct all $[0, 0, 1, 0, 0]$-admissible frequency matrices  $F$, together with the associated maxima matrices $M$,
\begin{equation}\label{E:(0,1,0)-admissible construction}
F=(\phi_0, \phi_1, \phi_2, \phi_3, \phi_4, \phi_5),
\qquad M=(\mu_0, \mu_1, \mu_2, \mu_3, \mu_4, \mu_5),
\end{equation}
with $\phi_0=(0, 0, 1, 0, 0)$, $\phi_1=(0, 0, 1, 0, \bullet)$, $\phi_2=(0, 0,  \bullet,  \bullet,  \bullet)$, $\phi_3=(\circ, \circ. \circ, \circ, \circ)$, \dots, and $\mu_i$ is $i$-th row in $M$. Denote by  $A(i)$  the set of all possible quintuples of frequencies in the $i$-th row, i.e. the set of quintuples of frequencies with the total sum $\leq k=1$. $A(0)$ is completely determined, and for the other we obviously have:
$$
\begin{aligned}
A(0)&=\{(0, 0, 1, 0, 0)\},\\
A(1)&=\{(0, 0, 1, 0, 0)\},\\
A(2)&=\{(0, 0, 0, 0, 0),  (0, 0, 1, 0, 0), (0, 0, 0, 1, 0), (0, 0, 0, 0, 1)\},\\
A(3)&=\{(0, 0, 0, 0, 0),  (1, 0, 0, 0, 0), (0, 1, 0, 0, 0), (0, 0, 1, 0, 0), (0, 0, 0, 1, 0), (0, 0, 0, 0, 1)\},\\
A(4)&=\{(0, 0, 0, 0, 0),  (1, 0, 0, 0, 0), (0, 1, 0, 0, 0), (0, 0, 1, 0, 0), (0, 0, 0, 1, 0), (0, 0, 0, 0, 1)\},\\
A(5)&=\{(0, 0, 0, 0, 0),  (1, 0, 0, 0, 0), (0, 1, 0, 0, 0), (0, 0, 1, 0, 0), (0, 0, 0, 1, 0), (0, 0, 0, 0, 1)\}.
\end{aligned}
$$
If we use printouts from the Python code 21AAIC.py in the Appendix, for $N=6$ and `highest weight' = [0, 0, 1, 0, 0], and if we activate ``print( 'i =', i, 'all fs =', all fs)'' on line 76, we get $A(1)$, \dots, $A(5)$ as above (with mixed-up left and right). We obtain the first six rows of $N_5^0$ if we activate ``print( 'i =', i, 'row1 =', row1)'' on line 73.

We may construct $F$ in steps $i=0,\dots, 4$, and by using (\ref{E:maxima matrix of F}) at each step we determine the rows of the associated maxima matrix
\begin{equation}\label{E:step (i+1) in the construction of F}
\phi_0, \phi_1, \dots,\phi_i; \
\mu_0, \mu_1, \dots,\mu_i  \rightsquigarrow
\phi_0, \phi_1, \dots,\phi_i, \phi_{i+1};  \
\mu_0, \mu_1, \dots,\mu_i, \mu_{i+1},
\end{equation}
$\phi_r\in A(r)$.
Since we are using (\ref{E:maxima matrix of F}) to determine $ \mu_{i+1}$, we have a map
$$
(\phi_{i+1},\mu_{i})\mapsto  \mu_{i+1},
$$
and if $ \mu_{i+1}$ does not satisfy the criteria (\ref{E:maxima matrix criteria}), we discard the newly constructed matrix from the further procedure. In Python code 21AAIC.py this function is defined on line 29 as ``filter frequencies(fs, ms1, *ks)''.

Since in our construction we need only the last maxima row $ \mu_{i}$ to determine whether the newly constructed frequency matrix in (\ref{E:step (i+1) in the construction of F})  is $[0, 0, 1, 0, 0]$-admissible, we could have discarded already used up $\mu_0, \mu_1, \dots, \mu_{i-1}$. On the other hand, if we are interested in the corresponding $[0, 0, 1, 0, 0]$-admissible partitions, we can keep track of the total contribution $\tau_i$ of rows $\phi_r$, $r\leq i$, to the corresponding  $[0, 0, 1, 0, 0]$-admissible partition,
$$
\tau_i=\sum_{r=1}^i\sum_{j=1}^w f_{rj}n_{rj}.
$$
Hence we should record our steps in the construction with data
\begin{equation}\label{E:step (i+1) in the construction of F 2}
\phi_0, \phi_1, \dots,\phi_i;\, \tau_i;\, \mu_i  \quad\rightsquigarrow\quad
\phi_0, \phi_1, \dots,\phi_i, \phi_{i+1};\, \tau_{i+1};\, \mu_{i+1}.
\end{equation}
Finally, if we are interested only in the number of constructed partitions, and do not need the constructed  $[0, 0, 1, 0, 0]$-admissible frequency matrices $\phi_0, \phi_1, \dots,\phi_i$,
we should record our steps in the construction only with data
\begin{equation}\label{E:step (i+1) in the construction of partitions}
 \tau_i;\, \mu_i  \quad\rightsquigarrow\quad
 \phi_{i+1};\, \tau_{i+1};\, \mu_{i+1},
\end{equation}
i.e. we should regard the step (\ref{E:step (i+1) in the construction of partitions}) as the map
\begin{equation}\label{E:step (i+1) in the construction of partitions 2}
(( \tau_i, \mu_i),  \phi_{i+1})\mapsto (\tau_{i+1}, \mu_{i+1}) ,
\end{equation}
and this is what the Python code 21AAIC.py does in two loops ``for total0, ms0 in all total ms0:'' and ``for fs1 in all fs:'' starting at lines 77 and 78. If in  21AAIC.py we change the ``While loop'' and the ``result'' into
\begin{verbatim}
    while True:
        all_total_ms1 = []
        row1 = get_row(i, w)
        all_fs = list(all_frequencies(i, *highest_weight))
        for total0, ms0 in all_total_ms0:
            for fs1 in all_fs:
                ms = filter_frequencies(fs1, ms0, *highest_weight)
                if ms is None:
                    continue
                total1 = row_fs_value(row1, fs1) + total0
                all_total_ms1.append((total1, ms))
        if i == N:
        #if i == w//2:
            break
        i += 1
        all_total_ms0 = all_total_ms1

result = []
for total1, ms1 in all_total_ms1:
    if total1 <= 6:
        result.append(total1)
print(result)
#print(len(all_total_ms1))
\end{verbatim}
then the activated ``print(len(all total ms1))'' on line 92 shows that there is altogether 164 final pairs $(\tau_5,\mu_5)$ in our construction of  $(0,1,0)$-admissible partitions with the support in the first $N=5$ rows. If we print only the list of $\tau_5\leq 6$ which appear in the final step of (\ref{E:step (i+1) in the construction of partitions}) we get
\begin{verbatim}
highest_weight = [0, 0, 1, 0, 0]
k = 1  w = 5
[0, 5, 6, 3, 4, 5, 6, 2, 6, 3, 4, 5, 1, 6, 5, 6, 4, 5, 6, 2, 6, 3]
\end{verbatim}
From this list we see that $\tau=6$ appears $7$ times, i.e. there is $7$ $(0,1,0)$-admissible partitions of $6$. The Python code 21AAIC.py is just a bit faster and more polished version of this code and gives for the number of all $(0,1,0)$-admissible partitions of $n\leq N=6$ the list
\begin{verbatim}
highest_weight = [0, 0, 1, 0, 0]
k = 1  w = 5
[[1, 1], [2, 2], [3, 3], [4, 3], [5, 5], [6, 7]]
\end{verbatim}
Of course, this example would have been much shorter if we simply wrote ``by hand'' all $(0,1,0)$-admissible partitions of $n\leq 6$, as we did for all $(2,0,0)$-admissible partitions of $n\leq 8$ in Example  \ref{Ex:(2,0,0)-admissible colored partitions}.
\end{example}

\begin{remark}\label{R:FFL theorem}
A basis of the finite dimensional representation $L_{C_\ell}(k_1, \dots, k_\ell)$ is parametrized in \cite{FFL} in terms of  symplectic Dyck paths---this result is closely related to our Conjecture
 \ref{C:conjecture}. For $\ell=4$ Feigin-Fourier-Littelmann's theorem can be rephrased in terms of
$(0,k_1,k_2,k_3,k_4))$-admissible $5\times9$ matrices\footnote{
Here $k_0=0$, but it works for any $k_0$ since  $L_{C_\ell}(k_1, \dots, k_\ell)\subset L_{C^{(1)}_\ell}(k_0, k_1, \dots, k_\ell)$.
}
$$
\begin{matrix}
k_0&0&k_1&0&k_2&0&k_3&0&k_4\\
k_0&0&k_1&0&k_2&0&k_3&0&\bullet\\
k_0&0&k_1&0&k_2&0&\bullet&\bullet&\bullet\\
k_0&0&k_1&0&\bullet&\bullet&\bullet&\bullet&\bullet\\
k_0&0&\bullet&\bullet&\bullet&\bullet&\bullet&\bullet&\bullet\\
\end{matrix}\,;
$$
in \cite{FFL} these matrices are written in the form
$$
\begin{matrix}
\bullet&\bullet&\bullet&\bullet&\bullet&\bullet&\bullet\\
&\bullet&\bullet&\bullet&\bullet&\bullet&\\
&&\bullet&\bullet&\bullet&&\\
&&&\bullet&&&\\
\end{matrix}.
$$

\end{remark}

\section{Appendix: Python code for counting admissible colored partitions}\label{Apendix}

Here we give the Python code which is in part explained in Example \ref{Ex: a construction of (0,1,0)-admissible partitions} where diagonals in the array (\ref{E:array of frequencies}) became the rows in the matrix (\ref{E:(0,1,0)-admissible construction}).

\begin{verbatim}

"""The program 21AAIC counts the number P(n) of [k1,k2,...,kw]
-admissible colored partitions of n <= N on the array Nw of w
rows of natural numbers. We input by hand N and k1,k2,...,kw as
the list 'highest_weight' on line 57. The result is a list of
pairs [n,P(n)]."""

def get_row(i, w):
    return [max(0, x) for x in range(i*2 - 1, i*2 - w - 1, -1)]

def all_subfrequencies(c, r):
    for g0 in range(r + 1):
        if c == 1:
            yield [g0]
        else:
            for gs in all_subfrequencies(c - 1, max(0, r - g0)):
                yield [g0] + gs

def all_frequencies(i, *ks):
    rest = []
    for x in reversed(ks):
        rest.append(x)
    if i > w//2:
        for gs in all_subfrequencies(w, k):
            yield gs
    else:
        for gs in all_subfrequencies(i*2-1, sum(ks[-i*2+1:])):
            yield gs + rest[2 * i-1:]

def filter_frequencies(fs, ms1, *ks):
    k = sum(ks)
    ms = []
    for j, f in enumerate(fs):
        if j:
            m = ms1[j - 1]
            m0 = ms[-1]
            if m0 > m:
                m = m0
            m += f
            if m > k:
                return None
            ms.append(m)
        else:
            ms.append(f)
    return ms

def row_fs_value(row, fs):
    s = 0
    for v, f in zip(row, fs):
        s += v * f
    return s

if __name__ == '__main__':
    """One should put by hand N and the 'highest wight'.
    For w=2n+1 for (k0,k1,...,kn) put [k0,0,k1,0,...,0,kn].
    For w=2n for (k0,k1,...,kn)e put [k0,k1,0,...,0,kn]."""
    N = 6
    highest_weight = [0, 0, 1, 0, 0]
    print("highest_weight =", highest_weight)
    w = len(highest_weight)
    k = sum(highest_weight)
    print('k =', k, ' w =', w)
    i = 1
    frequencies = {}
    result = []
    ms0 = []
    for j in range(0, len(highest_weight)):
        ms0.append(sum(highest_weight[-j-1:]))
    all_total_ms0 = [(0, ms0)]

    while True:
        all_total_ms1 = []
        row1 = get_row(i, w)
        #print( 'i =', i, 'row1 =', row1)
        min_next_row = get_row(i + 1, w)[-1]
        all_fs = list(all_frequencies(i, *highest_weight))
        #print( 'i =', i, 'all_fs =', all_fs)
        for total0, ms0 in all_total_ms0:
            for fs1 in all_fs:
                ms = filter_frequencies(fs1, ms0, *highest_weight)
                if ms is None:
                    continue
                total1 = row_fs_value(row1, fs1) + total0
                if total1 <= N:
                    if total1 > total0:
                        frequencies[total1] = 1 + frequencies.get(total1, 0)
                    if total1 <= N - min_next_row:
                        all_total_ms1.append((total1, ms))
        if row1[-2] > 0:
            result.append([row1[-1], frequencies.get(row1[-1], 0)])
            result.append([row1[-2], frequencies.get(row1[-2], 0)])
        if max(row1[-2:]) >= N:
            break
        i += 1
        all_total_ms0 = all_total_ms1

if w%2 == 0:
    result.remove([0,0])
print(result)
\end{verbatim}

\begin{example} In the Rogers-Ramanujan case, i.e.  for $w=2$ and $k=1$ we get:
\begin{verbatim}
highest_weight = [0, 1]

[1, 1], [2, 1], [3, 1], [4, 2], [5, 2], [6, 3], [7, 3], [8, 4],
[9, 5], [10, 6], [11, 7], [12, 9], [13, 10], [14, 12],
[15, 14], [16, 17], [17, 19], [18, 23], [19, 26], [20, 31];

highest_weight = [1, 0]

[1, 0], [2, 1], [3, 1], [4, 1], [5, 1], [6, 2], [7, 2], [8, 3],
[9, 3], [10, 4], [11, 4], [12, 6], [13, 6], [14, 8], [15, 9],
[16, 11], [17, 12], [18, 15], [19, 16], [20, 20].
\end{verbatim}
\end{example}

\begin{example} For $(2,1,0,0,1)^e$-admissible partitions in Example  \ref{Ex:some colored partitions for w even} we get:
\begin{verbatim}
highest_weight = [2, 1, 0, 0, 0, 0, 0, 1]

[1, 2], [2, 4], [3, 8], [4, 15], [5, 27], [6, 47], [7, 78],
[8, 128], [9, 205], [10, 323], [11, 499], [12, 763],
[13, 1148], [14, 1709], [15, 2516], [16, 3669],
[17, 5297], [18, 7589], [19, 10779], [20, 15204].
\end{verbatim}
\end{example}

\begin{example} Related to Remark \ref{R:FFL theorem}, we can modify slightly the above Python code for calculating dimensions of representations $L_{C_\ell}(k_1, \dots, k_\ell)$ (one way to do it is to use the change described in Example \ref{Ex: a construction of (0,1,0)-admissible partitions}, activate ``if i == w//2:'' on line 82 and ``print(len(all total ms1))'' on line 92, deactivate ``$\#i$f i == N:'' on line 81 and set highest weight = [0,0,1,0,1,0,2,0,2] for dim L [1, 1, 2, 2] ). For example, for $\ell=4$ we have:
\begin{verbatim}
dim L [1, 1, 2, 2]  =  3459456,
dim L [2, 1, 2, 2]  =  9848916,
dim L [0, 2, 2, 2]  =  4321512,
dim L [1, 2, 2, 2]  =  16358760,
dim L [2, 2, 2, 2]  =  43046721.
\end{verbatim}
\end{example}

\section*{Acknowledgement}
We thank Jim Lepowsky for numerous stimulating discussions and Shashank Kanade for writing the first program
(in Maple) for counting $(k,0,0)$-admissible colored partitions.

M.P. and A.M. would like to thank Universit\`a di Roma La Sapienza, il Dipartimento SBAI, for
the hospitality during our visit in September 2019, where this work originated.

M.P. is partially supported by Croatian Science Foundation under the project 8488 and by the QuantiXLie Centre of Excellence, a project cofinanced by the Croatian Government and European Union through the European Regional Development Fund - the Competitiveness and Cohesion Operational Programme (Grant KK.01.1.1.01.0004).

We thank one referee for remarks that led to an improvement of the exposition.


\end{document}